\documentclass[english,nolineno,a4paper]{socg-lipics-v2021}

\hideLIPIcs

\usepackage[dvipsnames]{xcolor}
\usepackage[utf8]{inputenc}

\usepackage{amsmath}
\usepackage{multicol}
\usepackage{graphicx}
\usepackage[colorinlistoftodos]{todonotes}

\usepackage{amsfonts}
\usepackage{enumitem}
\usepackage{tikz}
\usepackage{tikz-cd}
\usetikzlibrary{decorations.markings,positioning}
\usepackage{amsthm}
\usepackage{amssymb}
\usepackage{cite}
\usepackage{float}
\usepackage{microtype}
\usepackage{mathtools}
\usepackage{hyperref}
\usepackage{cleveref}
\usepackage{setspace}

\makeatletter
\def\blfootnote{\xdef\@thefnmark{*}\@footnotetext}
\makeatother

\DeclareMathSymbol{\minus}{\mathbin}{AMSa}{"39}

\title{Convergence of Leray Cosheaves for Decorated Mapper Graphs}

\authorrunning{Justin Curry, Washington Mio, Tom Needham, Osman Okutan, Florian Russold}

\author{Justin Curry}{University at Albany, State University of New York, USA }{jmcurry@albany.edu}{https://orcid.org/
0000-0003-2504-8388}{Supported by NSF CCF-1850052 and NASA 80GRC020C0016}

\author{Washington Mio}{Florida State University, Tallahassee, Florida}{wmio@fsu.edu}{
}{Supported by NSF grant DMS-1722995}

\author{Tom Needham}{Florida State University, Tallahassee, Florida}{tneedham@fsu.edu}{https://orcid.org/0000-0001-6165-3433}{Supported by NSF DMS-2107808}

\author{Osman Berat Okutan}{Max Planck Institute for Mathematics in the Sciences, Leipzig, Germany}{osman.okutan@mis.mpg.de}{
}{}

\author{Florian Russold\footnote{Corresponding Author}}{Graz University of Technology, Austria}{russold@tugraz.at}{}{Supported by the Austrian Science Fund (FWF): W1230}

\Copyright{Justin Curry, Washington Mio, Tom Needham, Osman Okutan, Florian Russold}



\funding{
 } 

\begin{document}

\maketitle

\begin{abstract}
We introduce decorated mapper graphs as a generalization of mapper graphs capable of capturing more topological information of a data set. A decorated mapper graph can be viewed as a discrete approximation of the cellular Leray cosheaf over the Reeb graph. We establish a theoretical foundation for this construction by showing that the cellular Leray cosheaf with respect to a sequence of covers converges to the actual Leray cosheaf as the resolution of the covers goes to zero.
\end{abstract}

\section{Introduction} \label{sec_intro}

Reeb\blfootnote{This is an abstract of a presentation given at CG:YRF 2023. It has been made public for the benefit of the community and should be considered a preprint rather than a formally reviewed paper. Thus, this work is expected to appear in a conference with formal proceedings and/or in a journal.} graphs and their discrete analogs---mapper graphs---are important tools in computational topology \cite{carr,carlsson,edelsbrunner,probabilisticconvergence}. 
They are used for data visualization \cite{breastcancer,extracting}, for comparing scalar fields via distances between Reeb graphs \cite{landi,bauer,categorifiedreebspaces}, and data skeletonization \cite{ge2011data}, among other things. 
Given a continuous map $f\colon X\rightarrow \mathbb{R}$, the Reeb graph $R_f$ summarizes the zero-dimensional connectedness of $X$. In practice, we deal with maps $f\colon P\rightarrow \mathbb{R}$ on discrete data sets $P$, where Reeb graphs are replaced by mapper graphs. A mapper graph is constructed by choosing a cover $\mathcal{U}=(U_i)_{i\in I}$ of $f(X)$ or $f(P)$, taking the components or clustering the data points of $f^{\minus 1}(U_i)$, collapsing each of the obtained components to a single vertex and connecting two vertices if their corresponding components have common points. In \cite{categorifiedreebspaces,convergence} it is shown that Reeb graphs can be viewed as cosheaves and that the cellular Reeb cosheaf w.r.t.\ a cover converges to the Reeb cosheaf if the resolution of the cover goes to zero, establishing a theoretical justification for working with finite covers. 
A limitation of this approach is that Reeb/mapper graphs fail to capture higher-dimensional topological features. To overcome these limitations we introduce decorated mapper graphs. 
\begin{figure}[t]
\centering
\includegraphics[width=12cm]{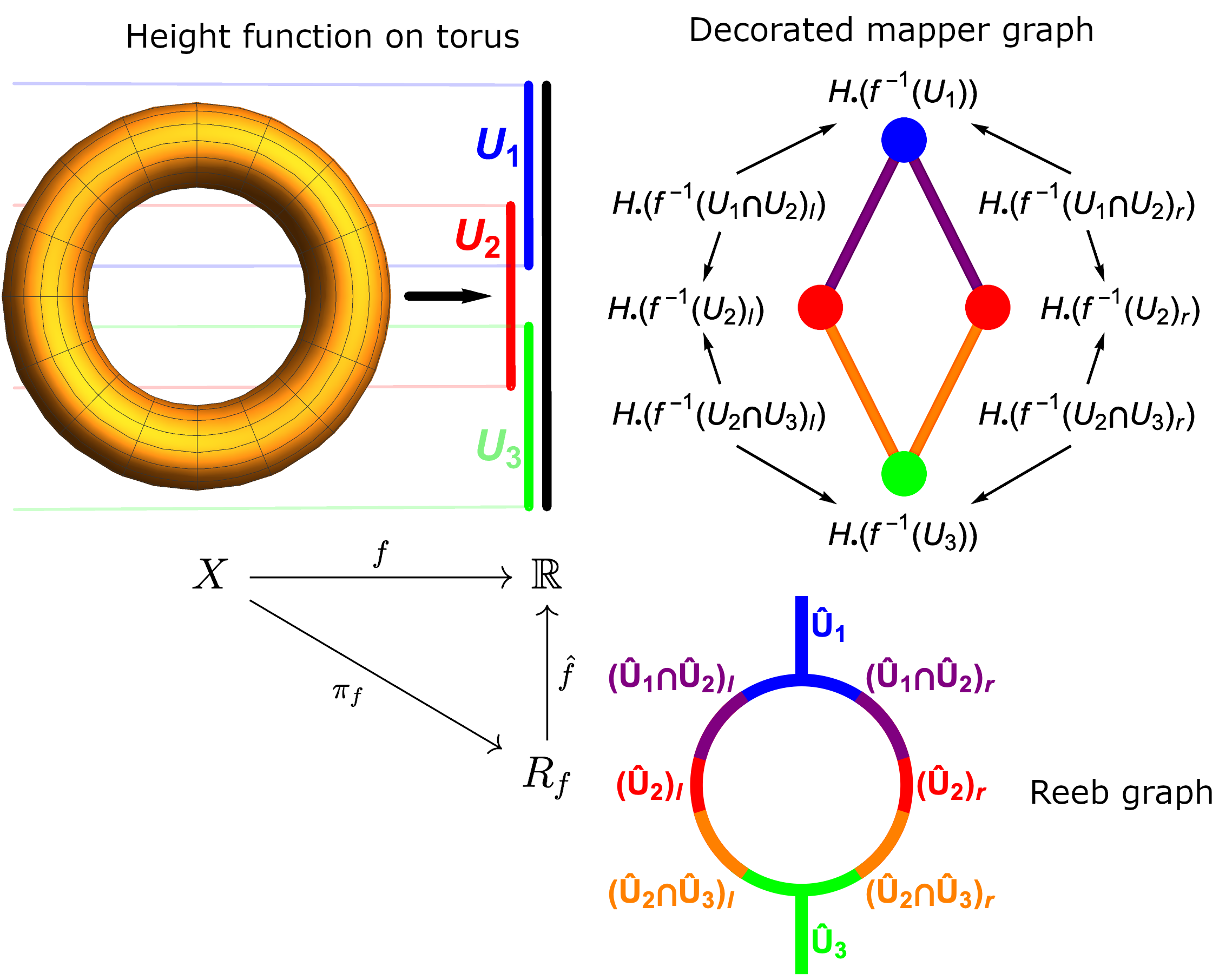}
\caption{Pulling back the cover $\mathcal{U}$ of $f(X)$ along the induced map $\hat{f}$ on the Reeb graph and refining it into components yields a cover $\hat{\mathcal{U}}$ of $R_f$. The decorated mapper graph is the nerve complex of this cover decorated by the homology $H_\bullet\big(\pi_f^{\minus 1}(-) \big)$ of the preimages of intersections in $\hat{\mathcal{U}}$ under the quotient map $\pi_f$. Since $f=\hat{f}\circ \pi_f$, this is analogous to using components of preimages under $f$.}
\label{fig_decorated_mapper}
\end{figure}
A decorated mapper graph collects the homology (possibly inferred from a point sample) in all degrees $H_\bullet(-)\coloneqq \bigoplus_{n\geq 0}H_n(-)$ of the components of $f^{\minus 1}(U_i)$ and $f^{\minus 1}(U_i\cap U_j)$ on the corresponding vertices and edges of the mapper graph; see Figure \ref{fig_decorated_mapper}. It can be viewed as a discrete approximation of the Leray cosheaf over the Reeb graph. 
In this abstract we show that the cellular Leray cosheaf w.r.t.\ a cover converges to the actual Leray cosheaf for any finite 1D cover $\mathcal{U}$ as the resolution of the cover goes to zero.   

\section{Leray Cosheaves} \label{sec_leray}

The Leray (pre)cosheaf \cite{curry,discretemorse,yoonghrist} parameterizes the homology of a space $X$ when viewed along a continuous map $f\colon X\rightarrow Y$. It does this by recording for each pair of open subsets $V\subseteq W \subseteq Y$ their homology and the induced map $H_\bullet\big(f^{\minus 1}(V)\big) \xrightarrow{H_\bullet(f^{\minus 1}(V)\subseteq f^{\minus 1}(W))} H_\bullet\big(f^{\minus 1}(W)\big)$.
We further assume that $X$ is compact and $H_\bullet\big(f^{\minus 1}(V)\big)$ is finite-dimensional for every $V\subseteq Y$.
\newpage
\begin{definition}[Leray (pre)cosheaf] \label{def_leray}
We define the graded Leray (pre)cosheaf $\mathcal{L}^f$ of a continuous map $f\colon X\rightarrow Y$ by the following assignments: For all $V\subseteq W\subseteq Y$ open
\begin{equation*} \label{eq_leray}
\begin{aligned}
\mathcal{L}^f(V)=\bigoplus_{n\in\mathbb{N}_0} \mathcal{L}^f_n(V)  &\coloneqq \bigoplus_{n\in\mathbb{N}_0}H_n\big(f^{\minus 1}(V)\big)\\
\mathcal{L}^f(V\subseteq W)=\bigoplus_{n\in\mathbb{N}_0} \mathcal{L}^f_n(V\subseteq W) &\coloneqq \bigoplus_{n\in\mathbb{N}_0}H_n\big(f^{\minus 1}(V)\subseteq f^{\minus 1}(W)\big) \hspace{2pt} .
\end{aligned}
\end{equation*}
\end{definition}
In practice, we can not access the whole Leray cosheaf. We can only get a cellular version \cite[4.1.6]{curry} given by its values on a finite cover of $f(X)$. An open cover $\mathcal{U}=(U_i)_{i\in I}$ of $f(X)$ defines a simplicial complex $\mathcal{N}_\mathcal{U}\coloneqq\{\sigma=(i_0,\ldots,i_k) \text{}|\text{ } \mathbf{U}_\sigma\coloneqq U_{i_0}\cap\ldots\cap U_{i_k}\neq\emptyset\}$, called the nerve complex of $\mathcal{U}$. We call $\mathcal{U}$ a \textbf{finite 1D cover} if it is finite and has a one-dimensional nerve complex and denote by $\leq$ the face-relation of $\mathcal{N}_\mathcal{U}$, i.e.\ $\sigma\leq\tau\iff\sigma\text{ is a face of }\tau$ .

\begin{definition}[Cellular Leray cosheaf] \label{def_discrete}
We define the cellular Leray cosheaf $D_\mathcal{U}\mathcal{L}^f$ w.r.t.\ a cover $\mathcal{U}$ as the cosheaf on $\mathcal{N}_\mathcal{U}$ given by the following assignments: For all $\sigma\leq \tau \in \mathcal{N}_\mathcal{U}$ 
\begin{equation*} \label{eq_discrete}
\begin{aligned}
D_\mathcal{U}\mathcal{L}^f(\sigma) &\coloneqq \mathcal{L}^f(U_\sigma) \\
D_\mathcal{U}\mathcal{L}^f(\sigma\leq\tau) &\coloneqq \mathcal{L}^f(U_\tau\subseteq U_\sigma) \hspace{2pt} .
\end{aligned}
\end{equation*}

\end{definition} 
If $\pi_f\colon X\rightarrow R_f$ is the quotient map from $X$ to the Reeb graph and $\hat{\mathcal{U}}$ a cover of $R_f$, then we define a decorated mapper graph as $D_{\hat{\mathcal{U}}} \mathcal{L}^{\pi_f}$; see Figure \ref{fig_decorated_mapper}.

\section{Convergence} \label{sec_convergence}

\begin{figure}[t]
\centering
\includegraphics[width=11cm]{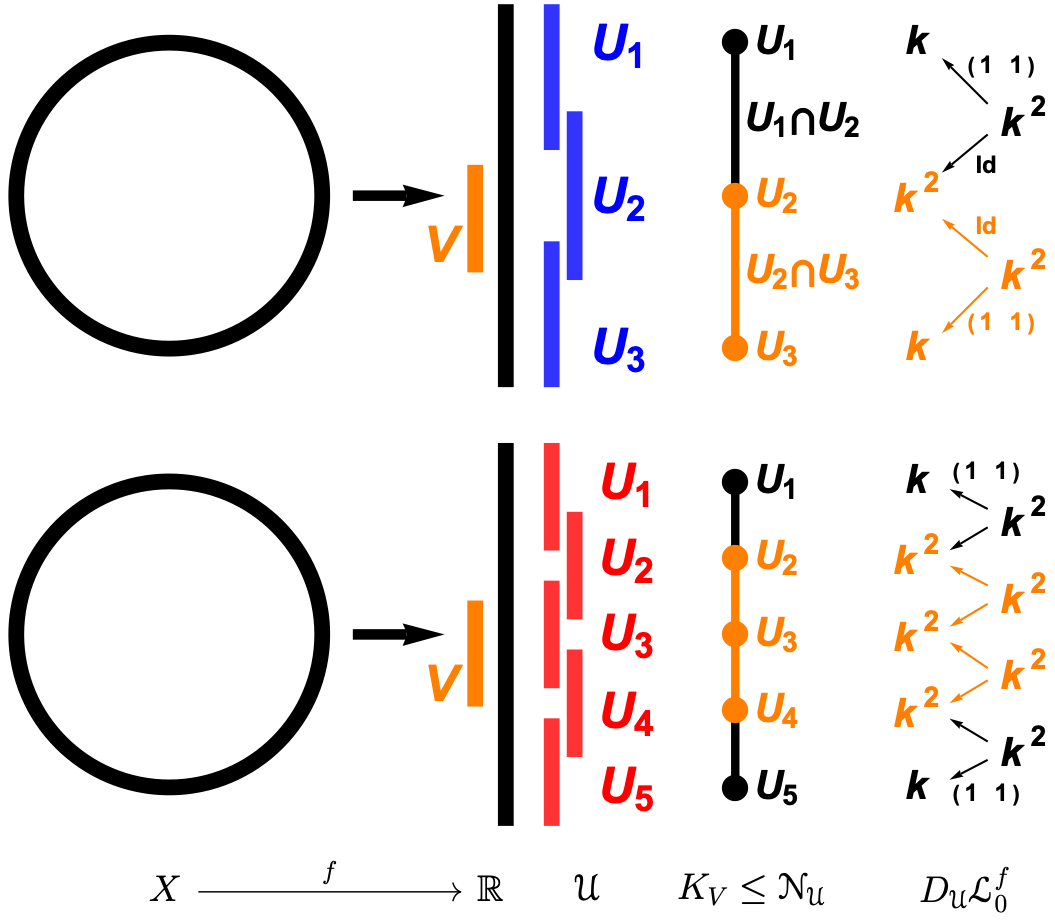}
\caption{The figure shows $D_\mathcal{U}\mathcal{L}^f_0\cong D_\mathcal{U}\mathcal{L}^f$ on $\mathcal{N}_\mathcal{U}$ w.r.t.\ two covers with different resolution as well as the restriction $D_\mathcal{U}\mathcal{L}^f_0|_{K_V}$ to $K_V$ (in orange). Homology is taken with coefficients in a field $k$ and unlabeled arrows represent identity maps. For the coarse cover we get $C_{\mathcal{U}}D_{\mathcal{U}}\mathcal{L}^f(V)\cong k\ncong H_0\big(f^{\minus 1}(V) \big)$ but for the finer one we get $C_{\mathcal{U}}D_{\mathcal{U}}\mathcal{L}^f(V)\cong k^2\cong H_0\big(f^{\minus 1}(V)\big)$.}
\label{fig_approximation_V}
\end{figure}

It is obvious that this process of discretization can lose information. This raises the question: How well is $\mathcal{L}^f$ represented by $D_\mathcal{U}\mathcal{L}^f$? 
To compare $\mathcal{L}^f$ and $D_\mathcal{U}\mathcal{L}^f$, we define a process of transforming $D_\mathcal{U}\mathcal{L}^f$ into a (pre)cosheaf on $Y$. We want to approximate the value of $\mathcal{L}^f$ on any open set $V\subseteq Y$ given only the information in $D_\mathcal{U}\mathcal{L}^f$. To this end, we define the subcomplex $K_V\coloneqq\{\sigma\in\mathcal{N}_\mathcal{U}\text{ }|\text{ } U_\sigma\cap V\neq\emptyset\}\leq\mathcal{N}_\mathcal{U}$, which can be viewed as a simplicial approximation of $V\subseteq \underset{\sigma\in K_V}{\bigcup}U_\sigma$; see Figure \ref{fig_approximation_V}. Moreover, if $V\subseteq W$, we obtain an inclusion $\iota\colon K_V \xhookrightarrow{} K_W$. 

\begin{definition}[Continuous extension] \label{def_continuous}
Let $\mathcal{U}$ be a finite 1D cover of $Y$. We define the continuous extension $C_\mathcal{U}D_\mathcal{U}\mathcal{L}^f$ of $D_\mathcal{U}\mathcal{L}^f$ as a (pre)cosheaf on $Y$ via the following assignments: For all $V\subseteq W\subseteq Y$ open
\begin{equation*} \label{eq_continuous}
\begin{aligned}
C_\mathcal{U}D_\mathcal{U}\mathcal{L}^f(V)&\coloneqq \bigoplus_{n\in\mathbb{N}_0}\Big(H_0\big(K_V;D_\mathcal{U}\mathcal{L}^f_n|_{K_V}\big)\oplus H_1\big(K_V;D_\mathcal{U}\mathcal{L}^f_{n\minus 1}|_{K_V}\big) \Big ) \\
C_\mathcal{U}D_\mathcal{U}\mathcal{L}^f(V\subseteq W)&\coloneqq\bigoplus_{n\in\mathbb{N}_0}\Big(H_0\big(\iota\big)_n\oplus H_1\big(\iota\big)_{n\minus 1} \Big )
\end{aligned}
\end{equation*}
where $H_i\big(K_V;D_\mathcal{U}\mathcal{L}^f_n|_{K_V}\big)$ is the degree $i$ homology \cite[6.2.2]{curry} of the cosheaf $D_\mathcal{U}\mathcal{L}^f_n$ restricted to $K_V$ and $H_i(\iota)_n$ is the map induced on cosheaf homology by $\iota\colon K_V \xhookrightarrow{} K_W$ cf.\ \cite[A.17]{russold2022persistent}. 
\end{definition}
As shown in \cite{curry,discretemorse,yoonghrist}, if $\mathcal{U}$ is a finite 1D cover of $f(X)$ the continuous extension yields $C_\mathcal{U}D_\mathcal{U}\mathcal{L}^f(V)\cong \underset{n\in\mathbb{N}_0}{\bigoplus}H_n\Big(f^{\minus 1}\big(\underset{\sigma\in K_V}{\bigcup}U_\sigma\big)\Big)$ approximating $H_\bullet\big(f^{\minus 1}(V) \big)$; see Figure \ref{fig_approximation_V}.
The following proposition shows that $C_\mathcal{U}D_\mathcal{U}\mathcal{L}^f(V\subseteq W)$ also gives us the correct induced 
map.

\begin{proposition} \label{thm_commutativity}
Let $\mathcal{U}$ be a finite 1D cover of $f(X)$. Then, for all $V\subseteq W\subseteq Y$ open, the following diagram commutes and the horizontal arrows are isomorphisms:
\begin{equation*} \label{eq_commutativity}
\begin{tikzcd}
H_0\big(K_V;D_\mathcal{U}\mathcal{L}_n^f|_{K_V}\big)\oplus H_1\big(K_V;D_\mathcal{U}\mathcal{L}_{n\minus 1}^f|_{K_V}\big) \arrow[r,"\cong"] \arrow[d,swap,"H_0(\iota)\oplus H_1(\iota)"] &[20pt] H_n\Big(f^{\minus 1}\big(\underset{\sigma\in K_V}{\bigcup}U_\sigma\big)\Big) \arrow[d,swap,"H_n\Big(f^{\minus 1}\big(\underset{\sigma\in K_V}{\bigcup}U_\sigma\big)\subseteq f^{\minus 1}\big(\underset{\sigma\in K_W}{\bigcup}U_\sigma\big)\Big)"] \\[30pt]
H_0\big(K_W;D_\mathcal{U}\mathcal{L}_n^f|_{K_W}\big)\oplus H_1\big(K_W;D_\mathcal{U}\mathcal{L}_{n\minus 1}^f|_{K_W}\big) \arrow[r,swap,"\cong"] & H_n\Big(f^{\minus 1}\big(\underset{\sigma\in K_W}{\bigcup}U_\sigma\big)\Big)
\end{tikzcd} \hspace{2pt} .
\end{equation*}
\end{proposition}
To talk about convergence, we have to define a distance on (pre)cosheaves. 
Assume $Y$ is a metric space and that all the involved (pre)cosheaves are constructible \cite[11.0.10]{curry}. 
For every $\epsilon>0$ define $V^\epsilon\coloneqq \{y\in Y \text{ }|\text{ } d(y,V)<\epsilon \}$. 
To a (pre)cosheaf $F$ on $Y$ we now associate $F^\bullet\coloneqq \{(F^\epsilon)_{\epsilon\geq 0}\hspace{2pt},\hspace{2pt}(F_\epsilon^\delta\colon F^\epsilon\rightarrow F^\delta)_{\epsilon\leq \delta}\}$, a parameterized family of (pre)cosheaves on $Y$, by setting $F^\epsilon(V)\coloneqq F(V^\epsilon)$. 
This allows us to define the distance of two (pre)cosheaves $F$ and $G$ as the interleaving distance of $F^\bullet$ and $G^\bullet$, i.e.\ $d(F,G)\coloneqq d_I(F^\bullet,G^\bullet)$. Denote by $\text{res}(\mathcal{U})\coloneqq \underset{U\in \mathcal{U}}{\text{sup}}\hspace{2pt}\underset{x,y\in U}{\text{sup}} d(x,y)$ the resolution of an open cover $\mathcal{U}$. We are now able to show that if the resolution of the cover $\mathcal{U}$ goes to zero, $C_\mathcal{U}D_\mathcal{U}\mathcal{L}^f$ converges to $\mathcal{L}^f$. 

\begin{theorem} \label{thm_interleaving}
If $f\colon X\rightarrow Y$ is continuous and $\mathcal{U}$ is a finite 1D cover of $f(X)$, then
\begin{equation*} \label{eq_interleaving}
d\big(C_\mathcal{U}D_\mathcal{U}\mathcal{L}^f,\mathcal{L}^f \big)\leq\text{res}(\mathcal{U}) \hspace{2pt} .
\end{equation*}
\end{theorem}
If we only consider $\mathcal{L}^f_0$, the degree-zero part of $\mathcal{L}^f$, the construction in Definition \ref{def_continuous} and Theorem \ref{thm_interleaving} specializes to an abelianization of the convergence result in \cite{convergence}.

\section{Future Work} \label{sec_future}

Theorem \ref{thm_interleaving} guarantees that, if we choose a cover of resolution $\leq\delta$, the cellular Leray cosheaf is a $\delta$-approximation of the continuous one. This result establishes a theoretical justification for working with finite covers. Since we have to deal with maps from finite point sets in practice, in future work we plan to investigate under what conditions we can infer the cellular Leray cosheaf of a map $f\colon X\rightarrow Y$ from a finite sample of $X$; cf.\ \cite{probabilisticconvergence}.

\bibliography{lib}

\end{document}